\documentclass[12pt]{amsart}
\usepackage{amssymb,latexsym,comment,url}

\usepackage{rotating}
\usepackage{rotating}
\usepackage{fancyhdr}
\usepackage{tikz}
\usepackage{tkz-euclide}
\usepackage{subcaption}
\usepackage[justification=centering]{caption}
\usetikzlibrary{calc}

\newcommand{\PP}{{\mathbb P}}
\newcommand{\Q}{{\mathbb Q}}

\newcommand{\rank}{\operatorname{rank}}

\newcommand{\Z}{{\mathbb Z}}

\newcommand{\Sage}{{\sf SAGE }}

\newenvironment{Proof}{\par\noindent{\sc Proof:}}%
                      {\hspace*{\fill}\nobreak$\Box$\par\medskip}
                       {\hspace*{\fill}\nobreak$\Box$\par\medskip}

\newtheorem{Proposition}{Proposition}[section]
\newtheorem{Theorem}[Proposition]{Theorem}
\newtheorem{Lemma}[Proposition]{Lemma}
\newtheorem{Corollary}[Proposition]{Corollary}

\theoremstyle{definition}

\newtheorem{Remark}[Proposition]{Remark}
 \newtheorem{Example}[Proposition]{Example}

\renewcommand{\baselinestretch}{1.1}

\begin{document}

\title[on rational triangles via algebraic curves]%
{on rational triangles via algebraic curves}

\author[M. Sadek]%
{Mohammad~Sadek}
\address{American University in Cairo, Mathematics and Actuarial Science Department, AUC Avenue, New Cairo, Egypt}
\email{mmsadek@aucegypt.edu}

\author[F. Shahata]%
{Farida~Shahata}
\address{American University in Cairo, Mathematics and Actuarial Science Department, AUC Avenue, New Cairo, Egypt}
\email{farida.mahmoud@aucegypt.edu}

\maketitle

\let\thefootnote\relax\footnote{Mathematics Subject Classification: 11G05, 11G30}
\begin{abstract}
A rational triangle is a triangle with rational side lengths. We consider three different families of rational triangles having a fixed side and whose vertices are rational points in the plane. We display a one-to-one correspondence between each family and the set of rational points of an algebraic curve. These algebraic curves are: a curve of genus 0, an elliptic curve, and a genus 3 curve. We study the set of rational points on each of these curves and describe some of its rational points explicitly.
\end{abstract}
\textbf{Keywords:} Rational triangles; Rational points; Algebraic curves.

\section{Introduction}

Several arithmetic questions on the geometry of the Euclidean plane have been a subject of interest in mathematical literature. One of the techniques to tackle such an arithmetic question is to construct a system of diophantine equations whose set of rational solutions provides an answer.
Therefore, answering these questions reveals an interplay between analytic geometry and arithmetic geometry.

 A rational subset of the plane is said to be {\em rational} if all pairwise distances between its points are rational. In 1945, Ulam wondered which rational subsets $S$ of the plane are infinite. One knows that a line and a circle have dense rational subsets. Erd\H{o}s conjectured that $S$ must be of a very special type. In \cite{Solymosi}, it was proved that the circle and the line are the only algebraic curves containing an infinite rational set. In order to prove the result, given an algebraic curve $C$ they constructed a new curve $C'$ of genus $g\ge 2$ whose number of rational points is larger than the size of any rational subset of the algebraic curve $C$. According to Faltings' Theorem, the number of rational points on $C'$ is finite, hence the result.

One may impose one more condition on rational sets and insist that their points are rational themselves, i.e., the $x$- and $y$- coordinates of the points are rational numbers. In \cite{Peeples}, it was proved that there exist infinitely many rational points on the $x$-axis lying at rational distances from 4 fixed points on the $y$-axis. The proof involves constructing an elliptic curve with positive rank. In \cite{Lorenzini}, it is shown that if the number of points on the $y$-axis is greater than 4, then there are only finitely many rational points on the $x$-axis lying at rational distances from the points on the $y$-axis. This holds because the number of such points is the number of rational points on an algebraic curve whose genus $g>1$.

A rational triangle is a triangle whose side lengths are rational numbers. One sees that the vertices of a rational triangle make up a rational set whose size is 3. An arithmetic question on rational triangles which has occupied a thrilling position since ancient times is which rational numbers appear as the area of a right-angled rational triangle. These rational numbers are called {\em congruent numbers}. A rational number $n$ is a congruent number if and only if the elliptic curve $y^2=x^3-n^2x$ is of positive rank.

A {\em Heron triangle} is a rational triangle whose area is rational. In \cite{Campbell,Duj, Fine, Luijk, Rusin}, Rational and Heron triangles with certain properties are investigated via studying algebraic curves and surfaces. Moreover, rational triangles are used to explore the size of the sets of rational points of some algebraic curves.

In this note, we are interested in rational triangles whose vertices are rational points in the plane. We consider three families of rational triangles with rational vertices.
The first family of such triangles consists of isosceles triangles. One knows that given a line segment $\ell$ of rational length and with rational endpoints, then there exist infinitely many isosceles triangles whose base is $\ell$, namely, those are the ones whose vertices are points on the perpendicular bisector of $\ell$. It turns out that out of the latter triangles, there are infinitely many whose sides are rational and the third vertex is a rational point in the plane. In other words, there are infinitely many rational points on the perpendicular bisector of $\ell$ which lie at a rational distance from the endpoints of $\ell$. This is proved by showing that these rational points on the perpendicular bisector are in one-to-one correspondence with the rational points on a genus 0 curve. In addition, we describe these points explicitly.

The second family of rational triangles consists of the triangles for which two of the vertices are the origin and a fixed rational point on the $x$-axis, whereas the third vertex is a rational point on a fixed line in the plane. These triangles possess rational areas, therefore, they are Heron triangles. We attach an elliptic curve of positive rank to this family, which implies the existence of infinitely many such triangles. We study this elliptic curve further by shedding some light on its torsion subgroup, and detecting conditions that force the rank to be at least 2.

The third family consists of the triangles for which two of the vertices are the origin and a fixed rational point $Q$ on the $x$-axis, and the third vertex is a rational point on the parabola $x=y^2$. It follows that such a triangle is a Heron triangle. Unlike the first and second families of triangles in this note, the third family turns out to consist of finitely many triangles. The reason is that the existence of a triangle in this family is equivalent to the existence of a rational point on a genus 3 curve. We display infinitely many rational points $Q$ for which the corresponding family of rational triangles is nonempty. The latter is achieved by exhibiting an explicit rational point on the corresponding genus 3 curve.

\section{Rational isosceles triangles and genus 0 curves}

We recall that a {\em rational triangle} is a triangle with rational side lengths. A point in the $xy$-plane is {\em rational} if its $x$- and $y$-coordinates are rational.

 Suppose we fix a rational point $P_1=(X_1,Y_1)$ in the $xy$-plane at a rational distance from the origin $O$ where $P_1$ and $O$ are distinct points. In this section we consider the problem of finding rational points $P=(X,Y)$ which lie at the same rational distance, say $R$, from the origin $O$ and the point $P_1$.

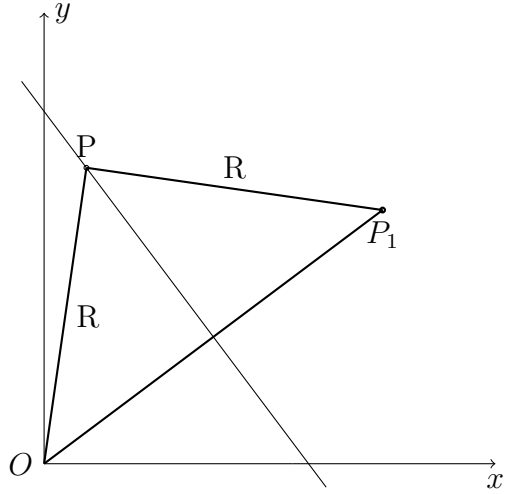
\begin{center}
\begin{tikzpicture}[scale=3]
\coordinate (O) at (0,0) coordinate[label = {left: $O$}] (20mm);
  \draw[-] (0,0) --  (1.1,0);
  \draw[->] (1.1,0) -- (2,0) coordinate[label = {below:$x$}] (20mm);

  \draw[->] (0,0) -- (0,2) coordinate[label = {right:$y$}] (20mm);

\draw (5/4,-5/48) -- (3/16, 21/16) circle [radius=0.3 pt]  node[anchor=south] {P};
\draw[black, thick] (0,0) -- (3/16,21/16) node[midway,right] {R};
\draw[black, thick] (1.5,9/8) -- (3/16,21/16) node[midway,above] {R};

\draw(3/16,21/16) -- (-1/10,407/240);
\draw [black, thick] (0,0) -- (1.5,9/8)  circle [radius=0.3 pt] coordinate[label = {below:$P_1$}] (20mm);

\end{tikzpicture}
\end{center}

In other words, we are seeking for the rational isosceles triangles $\triangle OPP_1$ with the fixed base $OP_1$ and the vertex $P$ itself being rational. We note that the point $P$ lies on the perpendicular bisector of the line segment joining $O$ and $P_1$, and that $\displaystyle R>\frac{\sqrt{X_1^2+Y_1^2}}{2} $. The rational point $P=(X,Y)$ together with the rational distance $R$ form a solution $(x:y:z:w)=(X:Y:R:1)$ of the following system of Diophantine equations

\begin{align}
\label{eq:intersection}
&x^2+y^2=z^2 \nonumber\\
&(x-X_1w)^2+(y-Y_1w)^2=z^2
\end{align}

The above equations describe the intersection curve of two quadratic surfaces in $\PP^3$. We will call this curve $E_{(X_1,Y_1)}$.
The curve $E_{(X_1,Y_1)}$ is birationally equivalent, hence isomorphic, to the singular cubic curve $C_{(X_1,Y_1)}$ defined by

\begin{equation*}
y^2=-16(X_1^2+Y_1^2)x(x+1)^2,
\end{equation*}
see for example Theorem 2.8 in \cite{Fsh}.

The curve $C_{(X_1,Y_1)}$ has a node at the point $S=(x,y)=(-1,0)$. There is an isomorphism of abelian groups between the nonsingular part $C_{(X_1,Y_1)}(\Q)\setminus \{S\}$ of $C_{(X_1,Y_1)}$ and the multiplicative group $\Q^{\times}$ of $\Q$, see \cite[Chapter 3, Exercise 3.5]{Silverman1}. In other words, one may find a rational parametrization for $C_{(X_1,Y_1)}\setminus\{S\}$.

Now we let $\mathcal R$ denote the set of all rational isosceles triangles $\triangle OPP_1$ with base $OP_1$ and the vertices $P$ and $P_1$ being rational where $|OP|=|PP_1|=R$. We note that a triangle in $\mathcal R$ is determined completely by the vertex $P$.

\begin{Proposition}
\label{prop1}
A triangle $\triangle OP_1P$ with base $OP_1$ lies in $\mathcal R$ if and only if there exists a $t\in\Q$ such that
\begin{eqnarray*}
X&=& \frac{X_1}{2} \pm \frac{t^2X_1(X_1^2+Y_1^2) - 4t(X_1^2 +Y_1^2)+4X_1}{2(t^2(X_1^2 +
Y_1^2) - 4)},\\
Y&=&\frac{Y_1}{2} \pm \frac{X_1(4t(X_1^2+Y_1^2)-t^2X_1(X_1^2+Y_1^2) - 4X_1)}{2Y_1(t^2(X_1^2+Y_1^2) -4)},\\
R&=&\frac{(t(X_1^2+Y_1^2) -2X_1)^2 +4Y_1^2}{2t^2Y_1(X_1^2+Y_1^2)-8Y_1}.
\end{eqnarray*}
\end{Proposition}
\begin{Proof}
The perpendicular bisector, $L$, of the line segment joining $O$ and $P_1=(X_1,Y_1)$ is described by
\begin{equation*}
 L:   y= \frac{-X_1}{ Y_1} x + \frac{X_1^2+Y_1^2}{2Y_1}
\end{equation*}
Since $P=(X,Y)$ lies on $L$ and the distance between $P$ and $O$ is $R$, i.e., $X^2+Y^2=R^2$, by substitution and solving a quadratic equation one gets
$$X= \frac{X_1}{2} \pm \frac{Y_1}{2} \sqrt{\frac{4R^2}{X_1^2+Y_1^2} -1} $$
Thus, for $X$ to be rational, we need to solve the following diophantine equation
\begin{equation}
\label{eq1}
\frac{4R^2}{X_1^2+Y_1^2} -1 = \delta ^2
\end{equation}


Since $X_1^2+Y_1^2 \neq 0$, we may assume without loss of generality that $Y_1\neq 0$. Now the point $\displaystyle (R,\delta)=\left(\frac{X_1^2+Y_1^2}{2Y_1},\frac{X_1}{Y_1}\right)$ is a rational solution for the latter diophantine equation. One concludes that (\ref{eq1}) has infinitely many rational solutions that can be parametrized as follows
\begin{eqnarray*}
(R,\delta)=\left( \frac{(t(X_1^2+Y_1^2) -2X_1)^2 +4Y_1^2}{2t^2Y_1(X_1^2+Y_1^2)-8Y_1},\frac{t^2X_1(X_1^2+Y_1^2) - 4t(X_1^2 +Y_1^2)+4X_1}{Y_1(4 -t^2(X_1^2 + Y_1^2)) }\right),\hspace{5mm} t\in\Q.
\end{eqnarray*}
Thus, $$X= \frac{X_1}{2} \pm \frac{t^2X_1(X_1^2+Y_1^2) - 4t(X_1^2 +Y_1^2)+4X_1}{2(t^2(X_1^2 +
Y_1^2) - 4)}.$$
\end{Proof}
 It is worth mentioning that for each $t\in\Q$, we obtain a distance $R$ and two values for $X$, which correspond to the points lying on the perpendicular bisector at a distance $R$ from the origin $ O$. So, the triangles that correspond to these two $X$-values are similar.
 \begin{Remark}
 In Proposition \ref{prop1}, we produced a parametric solution $(x:y:z:1)=(X:Y:R:1)$ to the curve of intersection of the two quadratic surfaces in (\ref{eq:intersection}).
 \end{Remark}
\begin{Corollary}
\label{cor1}
Given a rational point $P_1$ in the $xy$-plane lying at a rational distance from the origin $O$, there exists infinitely many rational points $P$ such that $\triangle OPP_1$ is a rational isosceles triangle for which $|OP|=|PP_1|$.
\end{Corollary}
\begin{Proof}
This is a direct consequence of Proposition \ref{prop1}.
\end{Proof}
\begin{Remark}
In corollary \ref{cor1}, one may replace the origin $O$ with any rational point $O'$ using a rational change of coordinates.
\end{Remark}
\begin{Example}
Taking $P_1=(3,4)$ and $t=0.5$, we get the point $\displaystyle P=\left(-\frac{25}{9}, \frac{125}{24}\right)$ which lies at a distance $R=\frac{425}{72}$ from both the origin and $P_1$. The triangle $\triangle OPP_1$ is rational.
\end{Example}

\section{Elliptic curves and rational triangles  }

In this section we investigate a family of rational triangles with a fixed side and link the number of these triangles with the size of the rational subgroup of an elliptic curve. We recall that a set $\mathcal S$ of points in the $xy$-plane is said to be a {\em rational distance set} if the distance between any two points in $\mathcal S$ is rational.

\subsection{The description of the problem}
Let $q,m,b$ be rational numbers. We consider the point $Q=(q,0)$ on the $x$-axis and the line $L:y=mx+b$. We are looking for the rational points $P=(X,Y)$ lying on $L$ at a rational distance, say $R$, from the origin $ O$; and at a rational distance, say $S$, from $Q$. We remark that this construction yields a rational distance set $\{O,P,Q\}$ of rational points.

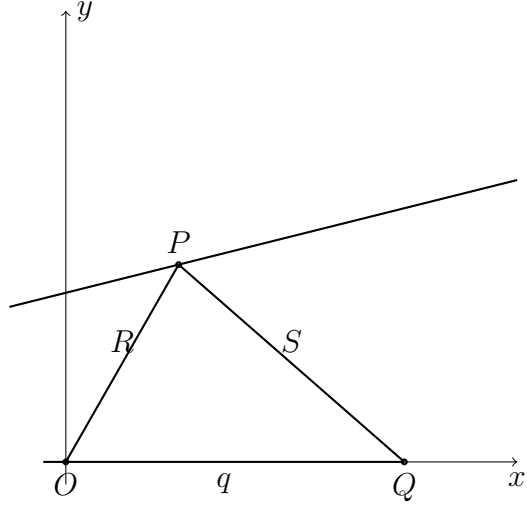
\begin{center}
\begin{tikzpicture}[scale=3]
\coordinate (O) at (0,0);
  \draw[black, thick] (-0.1,0) --  (1.5,0) circle [radius=0.3 pt] node[midway,below] {$q$} node[anchor=north] {$Q$};
  \draw[->] (1.5,0) -- (2,0) coordinate[label = {below:$x$}] (20mm);

  \draw[->] (0,-0.1) -- (0,2) coordinate[label = {right:$y$}] (20mm);

\draw [black, thick] (0,0) circle [radius=0.3 pt] node[anchor=north] {$O$}-- (1/2, 7/8) circle [radius=0.3 pt] node[midway,above] {$R$} node[anchor=south] {$P$};
\draw[black, thick] (1/2, 7/8) -- (1.5,0) node[midway,above] {$S$};

\draw [black, thick] (-0.25,11/16) -- (2,1.25);

\end{tikzpicture}
\end{center}

In other words, we are trying to find all the rational triangles $\triangle OQP$ with the fixed base $OQ$ and the vertex $P$ is lying on the line $L: y=mx+b$. We set $ A_{m,b,q}$ to be the set of all such triangles. We remark that any triangle in $A_{m,b,q}$ is a {\em Heron triangle}, i.e., it is a rational triangle with rational area. In fact, the area of a triangle $\triangle OQP$ is $|q(mX+b)|/2$.

The point $P=(X,Y)$ is a vertex of a triangle in $A_{m,b,q}$ if and only if $(x_1:x_2:x_3:x_4)=(X:R:S:1)$ is a solution of the following system of diophantine equations
\begin{eqnarray}
\label{eq:intersection2}
(1+m^2)x_1^2+2 mbx_1x_4+b^2x_4^2  =x_2^2 \nonumber\\
(1+m^2)x_1^2+2(mb-q)x_1x_4+(b^2+q^2)x_4^2=x_3^2
\end{eqnarray}
The above system represents an intersection curve $C_{m,b,q}$ of two quadratic surfaces in $\PP^3$. We note that, with a simple change of variables, $C_{0,b,q}$ has been extensively studied in relation with Heron triangles in \cite{Campbell}.

If $m\ne 0$, the point $\displaystyle \left(-\frac{b}{m}:\pm \frac{b}{m}: \pm \left( \frac{b}{m}  +q\right):1\right) $
is a rational point on the curve. By computing the discriminant of $C_{m,b,q}$, see \cite{Fsh}, we conclude that for the rational numbers $m,b$ and $q$ with $b \neq 0$, $q \neq 0$ and $q \neq \frac{-b}{m}$, the curve $C_{m,b,q}$ is an elliptic curve.

Let $C_{m,b,q}(\mathbb{Q})$ be the set of rational points on $ C_{m,b,q}$. We define the following relation on $C_{m,b,q}(\mathbb{Q})\setminus\{x_4=0\}$:
 The points $\left(\frac{x_1}{x_4}: \frac{x_2}{x_4}:\frac{x_3}{x_4}:1\right) \sim \left(\frac{x_1'}{x_4'}:\frac{x_2'}{x_4'}: \frac{x_3'}{x_4'}:1\right)$ if $\frac{x_1}{x_4}=\frac{x_1'}{x_4'}$, $\frac{x_2^2}{x_4^2}=\frac{x_2'^2}{x_4'^2}$ and $\frac{x_3^2}{x_4^2}=\frac{x_3'^2}{x_4'^2}$.
 This is clearly an equivalence relation and we may let $\mathcal C$ denote the set of equivalence classes.

\begin{Lemma}
\label{lem1}
 There exists a one-to-one correspondence between the set of triangles $A_{m,b,q}$ and the set of equivalence classes $\mathcal C$.
\end{Lemma}
\begin{Proof}
We define the one-to-one correspondence as follows. A triangle $\triangle OQP\in A_{m,b,q}$ will be sent to the equivalence class containing $(X:R:S:1)$ where $P=(X,Y)$, $R$ is the rational distance $|OP|$ and $S$ is the rational distance $|PQ|$. An equivalence class in $\mathcal C$, say represented by the rational point $(x_1:x_2:x_3:x_4)\in C_{m,b,q}(\Q), x_4\ne 0,$ will be sent to the triangle $OQP$ where $\displaystyle P=\left(\frac{x_1}{x_4},m\frac{x_1}{x_4}+b\right)$. This point lies at rational distance $\displaystyle \left|\frac{x_2}{x_4}\right|$ from $O$, and at rational distance $\displaystyle \left|\frac{x_3}{x_4}\right|$ from $Q$.
\end{Proof}

\begin{Remark}
The rational side lengths $q$, $R$ and $S$ of a triangle in $A_{m,b,q}$ force $\displaystyle X=\frac{q^2+R^2-S^2}{2q}$ to be rational itself.
Thus, $A_{m,b,q}$ is precisely the set of Heron triangles with fixed side of length $q$ and a vertex on the line $y=mx+b$.
\end{Remark}

 We know now that $C_{m,b,q}$ is an elliptic curve, except for finitely many possibilities for the values of $m$, $b$, and $q$, described as the intersection of two quadratic surfaces in $\PP^3$. We may wish to obtain a Weierstrass equation describing $C_{m,b,q}$. Indeed, when $b \neq 0$, $q \neq 0$ and $q \neq \frac{-b}{m}$, the curve $C_{m,b,q}$ is isomorphic to the elliptic curve $E_{m,b,q}$ described by the Weierstrass equation $y^2=x^3-27I_{m,b,q}x-27J_{m,b,q}$, see \cite{Fsh}, where

{\footnotesize\begin{eqnarray*}
I_{m,b,q}&=& 256(b^4 + 2b^3mq + (5m^2 + 4)b^2q^2 + 4mb(m^2+1)q^3 + (m^2 + 1)^2q^4),\\
J_{m,b,q}&=&4096(2b^2 + 2bmq + (m^2+1)q^2)(b^4 + 2b^3mq - b^2(7m^2 + 8)q^2 - 8bm(m^2+1)q^3  - 2(m^2+1)^2q^4).
\end{eqnarray*}}

 We recall that two elliptic curves defined over $\mathbb{Q}$ by the short Weierstrass equations $E_i: y^2=x^3+A_ix+B_i$ are isomorphic if and only if $\lambda^4 A_1=A_2$ and $\lambda^6 B_1=B_2$, for some nonzero $\lambda$ in $\mathbb{Q}$.

 We simply note that $I_{m,b,q}=b^4I_{m,1,\frac{q}{b}}$ and $J_{m,b,q}=b^6J_{m,1,\frac{q}{b}}$. So, given that $b \neq 0$ one obtains the following result.
  \begin{Proposition}
The elliptic curves $E_{m,b,q}$ and $\displaystyle E_{m,1,\frac{q}{b}}$ are isomorphic.
\end{Proposition}
  Thus, we will assume from now on that $b=1$.
Given $q$ and $m$ in  $\mathbb{Q}^*$ such that $q \neq \frac{-1}{m}$, we will write $C_{m,q}$ and $E_{m,q}$ instead of $C_{m,1,q}$ and $E_{m,1,q}$, respectively. The curve $C_{m,q}$ is given by the intersection of the following two quadratic surfaces, see (\ref{eq:intersection2}),

\begin{eqnarray}
\label{eq:intersection21}
(1+m^2)x_1^2 + 2mx_1x_4  + x_4^2&=& x_2^2\nonumber\\
(1+m^2)x_1^2 +2(m- q)x_1x_4  + (1+q^2)x_4^2 &=& x_3^2,
\end{eqnarray}
moreover, the elliptic curve $E_{m,q}$ is described by the following Weierstrass equation
{\footnotesize\begin{eqnarray}
\label{eq:Weierstrass}
y^2 = x^3&-&\frac{1}{3} \left[(m^2+1)^2 q^4 + 4m(m^2+1)q^3 + (5m^2+4)q^2 + 2m q + 1\right]x \nonumber\\&+&\frac{1}{27}\left[2(m^2+1)^2q^4 + 8m(m^2+1)q^3 + (7m^2+8)q^2 - 2m q - 1\right]\left[(m^2+1)q^2 + 2m q + 2\right].
\end{eqnarray}}
The explicit formulae for the isomorphism between the elliptic curves $C_{m,q}$ and $E_{m,q}$ are given in Appendix A.

\subsection{Torsion subgroup of $E_{m,q}(\Q)$}
We show that $E_{m,q}(\Q)$ contains the subgroup $\Z_2$. Moreover, we characterize those values for $m$ and $q$ such that $E_{m,q}(\Q)$ contains $\Z_2\times\Z_2$.

In order to find the points of order $2$ on $E_{m,q}$, we need to find the zeros of the polynomial $f_{m,q}(x)$ in the Weierstrass equation $y^2=f_{m,q}(x)$, (\ref{eq:Weierstrass}), describing the curve. Factoring the polynomial $f_{m,q}(x)$ yields that $\displaystyle x_1=\frac{1}{3}((m^2+1) q^2 + 2m q +2)$ is a zero of $f_{m,q}(x)$. In fact, this point $(x_1,0)$ of order $2$ corresponds to the rational point $\displaystyle\left(-\frac{1}{m}:  -\frac{1}{m}: -q-\frac{1}{m}:1\right)$ in $C_{m,q}(\Q)$. The existence of the point $(x_1,0)$ ensures that $E_{m,q}(\Q)$ contains $\Z_2$.

In fact, the other two zeros $x_2,x_3$ of $f_{m,q}(x)$ are given by
\[- \frac{1}{6}((m^2+1) q^2 + 2m q  +2) \pm  \frac{q}{2}\sqrt{(1+m^2)((1+m^2)q^2+4 m q + 4)}.\]
Therefore, that the torsion subgroup of $E_{m,q}(\Q)$ contains $\Z_2\times \Z_2$ is equivalent to $(1+m^2)((1+m^2)q^2+4 m q + 4)$ being a complete $\Q$-square.

For a fixed value of $m\in\Q$, one may parametrize the conic $t^2=(1+m^2)((1+m^2)q^2+4 m q + 4)$ and obtain the following value for $q$ which forces the torsion subgroup of $E_{m,q}(\Q)$ to contain $\Z_2\times\Z_2$
\[q:=q(n) = \frac{1 - m^2(2+3m^2) + 2 (1+m^2) n  + n^2}{m (1 + m^2(2+m^2) - n^2)},\hspace{5mm}n\in\Q.\]
\begin{Proposition}
The torsion subgroup of the elliptic curve $E_{m,q}(\Q)$ contains the cyclic group $\Z_2$. Furthermore, for a fixed value of $m$, there exist infinitely many values for $q$ such that the torsion subgroup of $E_{m,q}(\Q)$ contains $\Z_2\times\Z_2$.
\end{Proposition}

Now we will introduce infinitely many values for $q$ such that the torsion subgroup of $E_{m,q}$ contains a point of order $4$.
\begin{Proposition}
\label{prop:pointorder4}
For some nonzero rational $m$, if $\displaystyle q=\frac{4(2t-m)}{1+m^2-4t^2}$ for some $t\in\Q$, such that $\displaystyle t \neq m \pm \frac{\sqrt{1+m^2}}{2}$, then the point $P=(x,y)$ given by
\begin{eqnarray*}
x&=&(m^4 - 24m^3t + (104t^2 + 6)m^2 -mt(160 t^2+24) + 80t^4 + 24t^2 + 5)/3(m^2 - 4t^2 + 1)^2\\
y&=& -2(3m^2 - 8mt + 4t^2 - 1)(m^2 - 4mt + 4t^2 + 1)/(m^2 - 4t^2 + 1)^2
\end{eqnarray*}
is a point of order 4 in $E_{m,q}(\mathbb{Q})$.
 \end{Proposition}
 \begin{Proof}
  We first note that if $\displaystyle t = m \pm \frac{\sqrt{1+m^2}}{2}$ then $q=-\frac{1}{m}$ and hence $E_{m,q}$ is a singular curve.
 The order of $P$ can be easily checked using the duplication formula on $E_{m,q}$, see for example \cite[p. 59]{Silverman1}. In fact, $\displaystyle x(2P)=  \frac{2(4t^4 - 8t^3 + 16t^2 - 12t + 3)}{3(2t^2 - 1)^2}$ and $y(2P)=0$.
 \end{Proof}
 \begin{Remark}
The point $P$ in Proposition \ref{prop:pointorder4} is in correspondence with the following rational point in $C_{m,q}(\Q)$, see (\ref{eq:intersection21}),
{\footnotesize$$\left(\frac{q}{2}:\sqrt{\frac{1+m^2}{4} q^2+ m q+1}:\sqrt{\frac{1+m^2}{4} q^2+ m q+1}:1\right)$$
$$=\left(2\frac{2t-m}{1+m^2-4t^2}:\frac{m^2 - 4mt + 4t^2 + 1}{m^2 - 4t^2 + 1}:\frac{m^2 - 4mt + 4t^2 + 1}{m^2 - 4t^2 + 1}:1\right).$$}
The latter point corresponds to the isosceles triangle in the set $A_{m,q}$ whose base is $|q|$, and the two other sides are of length $\displaystyle\frac{m^2 - 4mt + 4t^2 + 1}{m^2 - 4t^2 + 1}$. Note that the area of this triangle is $\displaystyle \left|\frac{q(q m+2)}{4}\right|\in\Q$.
 \end{Remark}

\subsection{Rank of $E_{m,q}$}
We are going to investigate the rank of $E_{m,q}(\Q)$, and see the impact of the positivity of the rank on the size of the set $A_{m,q}$.

\begin{Theorem}
\label{thm1}
One has $\rank E_{m,q}(\mathbb{Q}(m, q)) \ge 1$.
\end{Theorem}
\begin{Proof}
The point $\displaystyle P_{m,q}=\left(-\frac{(2m^2-1)q^2 + 4 m q + 1}{3},  q(mq + 1)^2\right)$ is a point in $E_{m,q}(\mathbb{Q}(m, q))$.
By specializing with the values $q=1$ and $m=1$, the elliptic curve $E_{1,1}(\mathbb{Q})$ is defined by $y^2=x^3-8x+8$ where $P_{1,1}=(-2,4)\in E_{1,1}(\Q)$ is of infinite order. Therefore, by Silverman's Specialization Theorem, \cite[Appendix C, \S 20]{Silverman1}, $\rank E_{m,q}(\mathbb{Q}(m,q))\ge 1$.
\end{Proof}

\begin{Corollary}
For all but finitely many pairs $(m,q)\in\Q\times\Q$, the set $A_{m,q}$ contains infinitely many rational triangles with base $|q|$ and whose vertices are lying on the line $y=mx+1$.
\end{Corollary}
In what follows we choose $m$ and $q$ so that we can construct elliptic curves $E_{m,q}$ with higher rank.

\begin{Theorem}
\label{thm2}
Set $q:=q(h)=\frac{1-h^2}{2h}$. One has $\rank E_{m,q}(\Q(m,h))\ge 2$.
\end{Theorem}
\begin{Proof}
The point $Q_{m,q}=(x_{m,h},y_{m,h})$ given by
\begin{eqnarray*}
x_{m,h}&=&\frac{h^4(1+m^2) - 4h^3m - 2h^2(m^2+3m -3)  + 4h(m+3)  + m^2 + 6m+1}{12h^2}\\
 y_{m,h}&=&\frac{(m h^2 - 2h - m)(h + 1)(hm - m - h - 1)}{4 h^3}
 \end{eqnarray*}
 lies in $E_{m,q}(\Q(m,h))$. Moreover, the point $P_{m,q}$, see Theorem \ref{thm1}, is a point in $E_{m,q}(\Q(m,h))$. Taking $m=1$ and $h=1/2$, one has $P_{1, 3/4}=(-73/48 ,147/64)$  and $Q_{1, 3/4}=(121/24, 21/2)$. Further, using \Sage one obtains that these points are linearly independent in $E_{1,3/4}$. The result now follows using Silverman's Specialization Theorem.
\end{Proof}
\begin{Remark}
\label{rem:congruent}
The point $Q_{m,q(h)}\in E_{m,q(h)}$ in Theorem \ref{thm2} corresponds the point $(0:1: \frac{1+h^2}{2h}:1)\in C_{m,q(h)}.$ The triangle in $A_{m,q(h)}$ that corresponds to $Q_{m,q(h)}$ is a right Heron triangle with area: $\displaystyle A(h)=\left|\frac{1-h^2}{4h}\right|$. Thus, the numbers $A(h)$ are congruent numbers.
\end{Remark}
\begin{Theorem}
\label{thm3}
Let $\displaystyle q:=q(u,m)=\frac{2(u-m)}{1+m^2-u^2}$. Then $\rank E_{m,q}(\Q(m,u))\ge 2$.
\end{Theorem}
\begin{Proof}
The point $H_{m,q}=(x_{m,u},y_{m,u})$ given by
\begin{eqnarray*}
x_{m,u}&=&\frac{5m^4 +(6- 10u)m^3 + (4u^2 -6u+10)m^2  + 2(u^3 -3u^2-5u+3)m - u^4 + 6u^3 - 6u + 5}{3(m^2 - u^2 + 1)^2}\\
y_{m,u}&=&\frac{2(m - u + 1)^2(m^2 - um + 1)}{(m^2 - u^2 + 1)^2}
\end{eqnarray*}
together with the point $P_{m,q}$ are two linearly independent points. This can be proved by specializing $m=1,u=3$ and hence $q=-4/7$.
\end{Proof}
\begin{Remark}
The point $H_{m,q}$ in Theorem \ref{thm3} in $E_{m,q}(\Q(m,u))$ is in correspondence with the point
\[\left(\frac{2(u-m)}{1+m^2-u^2}: \frac{m^2 - 2mu + u^2 + 1}{m^2 - u^2 + 1}: -\frac{(m-1-u)(m+1-u)}{m^2 - u^2 + 1}:1 \right)\]
on $C_{m,q}(\Q(u,m))$. The latter point gives rise to a right rational triangle whose base is $|q|$ and area is $\displaystyle \left|\frac{(m - u)(m - u - 1)(m - u + 1)}{(m^2 - u^2 + 1)^2}\right|$, in particular this area is a congruent number.
\end{Remark}

\section{Rational triangles via curves of genus 3}
The problem we discuss in this section is similar to the one in \S 4. We let $q$ be a rational number. We consider the point $Q=(q,0)$ on the $x$-axis and the parabola $x=y^2$. We are searching the plane for the rational points $P=(X,Y)$ lying on the parabola at a rational distance, say $R$, from the origin $ O$; and at a rational distance, say $S$, from $Q$. Again this construction yields a rational distance set $\{O,P,Q\}$ of rational points.

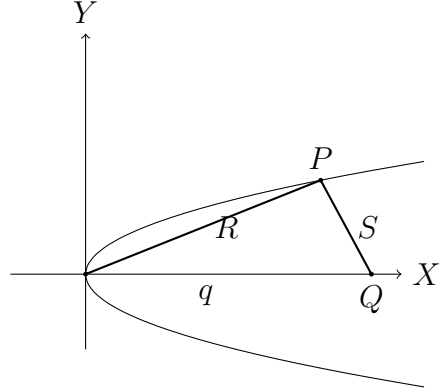
\begin{center}
\begin{tikzpicture}
  \draw[->] (-1,0) -- (4.2,0) node[right] {$X$}  node[midway,below] {$q$} ;
  \draw[->] (0,-1) -- (0,3.2) node[above] {$Y$};

  \draw[scale=0.5,domain=-3:3,smooth,variable=\y]  plot ({\y*\y},{\y});
  \draw[black, thick](0,0) circle [radius=0.5 pt] node[anchor=north] {} -- (25/8, 5/4) circle [radius=0.5 pt] node[midway,right] {$R$} node[anchor=south] {$P$};
  \draw[black, thick](25/8, 5/4) -- (3.8, 0) circle [radius=0.5 pt] node[midway,right] {$S$} node[anchor=north] {$Q$};

\end{tikzpicture}
\end{center}

The set $S_q$ is the set of rational triangles whose base is $OQ$ and whose third vertex $P$ is lying on the parabola $x=y^2$. We remark that $\triangle OQP$ is a Heron triangle. Furthermore, a point $P=(X,Y)$ is a vertex of some triangle in $S_q$ if and only if $(x_1:x_2:x_3:x_4:x_5)=(X:Y:R:S:1)$ is a rational solution for the following intersection $C_q$ of quadratic surfaces in $\PP^4$
\begin{eqnarray*}
x_1^2+x_2^2 =& x_3^2 \\
(x_1-qx_5)^2+x_2^2 =&  x_4^2 \\
x_2^2=& x_1x_5\\
\end{eqnarray*}
where $R$ is the distance between $P$ and $O$, and $S$ is the distance between $P$ and $Q$. The point $(0:0:0:\pm q:1)$ lies in $C_q(\Q)$. Yet, this point corresponds to a degenerate triangle in $S_q$.

It is known that the genus of a smooth complete intersection of three quadratic surfaces in $\PP^4$ is $5$. However, using the Jacobian criterion of smoothness, the intersection curve $C_q$ has two ordinary double points, namely $(0:0:0:\pm q:1)$. It follows that the curve $C_q$ is of genus 3. Faltings' celebrated theorem implies that the set of rational points $C_q(\Q)$ is finite as the genus of $C_q$ is greater than $1$.
This yields the following result.
\begin{Corollary}
The set $S_q$ is finite.
\end{Corollary}
One may ask whether the set $S_q$ can be nonempty.
\begin{Theorem}
\label{thm4}
If $\displaystyle q=\frac{(u^2+1)^2}{8u^2}$ for some $u\in\Q\setminus\{0,\pm1\}$, then $S_q$ contains an isosceles Heron triangle. Similarly, if $\displaystyle q=\frac{(u^2-1)^2}{4u^2}$ for some $u\in\Q\setminus\{0,\pm1\}$, then $S_q$ contains a right Heron triangle.
\end{Theorem}
\begin{Proof}
When $q=(u^2+1)^2/8u^2$, one sees that \[ (x_1:x_2:x_3:x_4:x_5)=\left(\frac{(u^2-1)^2}{4u^2}:\frac{u^2-1}{2u}:\frac{u^4 - 1}{4 u^2}: \frac{(u^2 + 1)^2}{8u^2}:1\right)\]
is a rational point in $C_q(\Q)$. This produces an isosceles rational triangle in $S_q$ as $q=S$.

If $q=(u^2-1)^2/4u^2$, then \[\left(\frac{(u^2-1)^2}{4u^2}:  \frac{u^2-1}{2u}:   \frac{u^4 - 1}{4 u^2}:  \frac{u^2-1}{2u}:1\right)\]
is a rational point in $C_q(\Q)$. Therefore, one has a right Heron triangle in $S_q$ as $q^2+S^2=R^2$.
\end{Proof}
\begin{Remark}
In Theorem \ref{thm4}, we produced the congruent numbers $\displaystyle\frac{(u^2 - 1)^3}{16u^3}$ which is the area $q S$ of the right Heron triangle.
\end{Remark}

\appendix
 \section{The isomorphism between $E_{m,q}$ and $C_{m,q}$}

 Recall that the curve $C_{m,q}$ is defined by
 \begin{eqnarray*}
(1+m^2)x_1^2 + 2mx_1x_4  + x_4^2= x_2^2,
 (1+m^2)x_1^2 +2(m - q)x_1x_4  + (1+q^2)x_4^2 = x_3^2.
\end{eqnarray*}
The curve $C_{m,q}$ is isomorphic to the curve $C'_{m,q}$ defined by
\begin{equation*}
y^2=\frac{(1+m q)^2}{4} x^4+ +q(1+m q)x^3z +((1-\frac{m^2}{2})q^2-m q+\frac{1}{2})x^2 z^2-q(1+m q)x z^3 +\frac{(1+m q)^2}{4} z^4
\end{equation*}
via the following isomorphism $\psi:C_{m,q}\to C'_{m,q}$; $(x_1:x_2:x_3:x_4)\mapsto (x:y:z)$ defined by
{\footnotesize\[
\psi(x_1:x_2:x_3:x_4) =
  \begin{cases}
      (m x_1+x_4:  x_3(x_1+x_2) : x_1+x_2 )   & \text{ if }(x_1:x_2:x_3:x_4) \neq \left(\frac{-1}{m}:\frac{1}{m}:q+\frac{1}{m}:1\right)\\
      (1: \frac{m q +1}{2}: 0) & \text{ if }(x_1:x_2:x_3:x_4)=\left(\frac{-1}{m}:\frac{1}{m}:q+\frac{1}{m}:1\right)
  \end{cases}
\]} whereas the inverse is given by

{\footnotesize\[
\psi^{-1}(x:y:z) =
  \begin{cases}
      \left(\frac{z^2-x^2}{2z}: \frac{x^2 + z^2}{2z}: \frac{y}{z} : \frac{m(x^2-z^2)+ 2x z}{2z}\right) & \text{ if }(x:y:z) \neq \left(1:\pm \frac{mq+1}{2}: 0 \right)\\
      \left(\frac{-1}{m}:\frac{1}{m}:q+\frac{1}{m}:1\right) & \text{ if } (x:y:z) =\left(1:\frac{mq+1}{2}: 0 \right) \\
       \left(\frac{-1}{m}:\frac{1}{m}:-q-\frac{1}{m}:1\right) & \text{ if } (x:y:z) =\left(1:\frac{-mq-1}{2}: 0 \right).
  \end{cases}
\]
 }

 Now the curve $C'_{m,q}$ is isomorphic to the elliptic curve $E_{m,q}$ defined by the following Weierstrass equation
 {\footnotesize\begin{eqnarray}
Y^2 = X^3&-&\frac{1}{3} \left[(m^2+1)^2 q^4 + 4m(m^2+1)q^3 + (5m^2+4)q^2 + 2m q + 1\right]X \nonumber\\&+&\frac{1}{27}\left[2(m^2+1)^2q^4 + 8m(m^2+1)q^3 + (7m^2+8)q^2 - 2m q - 1\right]\left[(m^2+1)q^2 + 2m q + 2\right]\nonumber.
\end{eqnarray}}
Let $\displaystyle G_{m,q}=\{((m^2q^2+2mq+q^2 -1)/3  ,\pm q)\}\subset E_{m,q}(\Q)$, and $\displaystyle G'_{m,q}=\{(1: \pm \frac{mq+1}{2}:0)\}\subset C'_{m,q}(\Q)$.
The isomorphism $\phi:C'_{m,q}\to E_{m,q}$ is defined as follows: $\phi(G'_{m,q})=G_{m,q}$; otherwise $\phi(x:y:1) =(X, Y)$ where
{\footnotesize\begin{eqnarray*}
X&=&y(1+m q)+\frac{1}{2} (m^2q^2 +m q+1  )x^2+q(m q+1)x+\frac{1- 2 m q + 2q^2 - m^2q^2}{6}\\
Y&=&\frac{m q + 1}{2}(  2q y + 2(1+m q)x y + (m q+1)^2x^3  + 3q(1+m q)x^2 + x(1-m^2q^2+2q^2-2m q) - q( m q+1)).
\end{eqnarray*}
}
Whereas the inverse is defined by $\phi^{-1}(G_{m,q})=G'_{m,q}$; otherwise $\phi^{-1}(X,Y)=(x:y:1)$ where
\begin{eqnarray*}
x& =& -\frac{3Y +( 2- 3X)q + 2m q^2 + (1+m^2)q^3}{ (m q + 1)(b(m,q) - 3X )}\\
y&=&\frac{54X^3 -27 b(m,q) X^2 - 27Y^2 - 54q Y + a(m,q)}{6(m q + 1)( b(m,q) -3X)^2}
\end{eqnarray*}
where
\footnotesize{\begin{eqnarray*}
a(m,q)&=& q^6(m^2 + 1)^3 + 6m q^5(m^2 + 1)^2+3q^4(3m^2 - 1)(m^2 + 1)-4m q^3(m^2 + 3)-3q^2(3m^2 + 8)+ 6m q-1,\\
b(m,q)&=&m^2q^2 + 2m q + q^2 - 1.
\end{eqnarray*}}

\end{document}